\documentclass[11pt]{amsart}
\usepackage{verbatim}
\usepackage{bm} 
\usepackage{amssymb}
\usepackage{colonequals}
\usepackage{mathtools}
\usepackage[bookmarks,colorlinks,breaklinks]{hyperref}  
\usepackage[top=2cm, bottom=4.5cm, left=3cm, right=3cm]{geometry}
\hypersetup{linkcolor=blue,citecolor=blue,filecolor=blue,urlcolor=blue} 
\usepackage{color}
\usepackage{amsfonts,amsmath,amsthm,amssymb,amscd, hyperref}

\numberwithin{equation}{section}

\newtheorem{theorem}{Theorem}[subsection]
\renewcommand{\thetheorem}{%
  \ifnum\value{subsection}=0 \thesection.\arabic{theorem}%
  \else \thesubsection.\arabic{theorem}%
  \fi}
\newtheorem{corollary}[theorem]{Corollary}
\newtheorem{lemma}[theorem]{Lemma}
\newtheorem{proposition}[theorem]{Proposition}

\theoremstyle{definition}

\newtheorem{definition-theorem}[theorem]{Definition-Theorem}

\newtheorem{remark}[theorem]{Remark}

\theoremstyle{remark}

\newtheorem*{remark*}{Remark}

\usepackage{enumitem}
\setenumerate[1]{label=\textup{(\alph*)}}
\setenumerate[2]{label=\textup{(\roman*)}}

\renewcommand{\l}{\lambda}

\newcommand{\bP}{\mathbb{P}}

\newcommand{\id}{{\rm id}}
\newcommand{\into}{\hookrightarrow}

\newcommand{\N}{{\mathbb N}}

\newcommand{\A}{{\mathbb A}}
\newcommand{\bA}{{\mathbb A}}

\newcommand{\Qbar}{{\overline{\mathbb Q}}}

\newcommand{\Kbar}{{\overline{K}}}

\newcommand{\hhat}{\widehat{h}} 
\newcommand{\lra}{\longrightarrow}

\newcommand{\bt}{{\mathbf t}}
\newcommand{\bl}{{\bm{\lambda}}}
\newcommand{\PP}{{\mathbb P}}

\newcommand{\cL}{\mathcal{L}}

\begin{document}
\raggedbottom
\title[Collision parameters]{Collision Parameters in Higher-Dimensional Families of Polynomials}

\author{Dragos Ghioca}
\address{Department of Mathematics, University of British Columbia, Vancouver, BC V6T 1Z2}
\email{dghioca@math.ubc.ca}

\author{Negin Shadgar}
\address{Department of Mathematics, University of British Columbia, Vancouver, BC V6T 1Z2}
\email{negin@math.ubc.ca}



\begin{abstract}
Let $d\geq 3$ and consider the generic family of monic normalized polynomials
\[
  f_{\bt}(z)=z^d+t_1z^{d-2}+t_2z^{d-3}+\cdots+t_{d-2}z+t_{d-1},
  \qquad \bt:=(t_1,\dots, t_{d-1}).
\]
Given pairwise distinct starting points $\alpha_1,\ldots,\alpha_d\in\Qbar$ and a target point $\beta\in\Qbar$, we show that the set of parameters $\bl\in\A^{d-1}(\Qbar)$ for which $\beta$ belongs to the forward orbit of every $\alpha_i$ under $f_{\bl}$ is not Zariski dense in $\A^{d-1}$. 
\end{abstract}

\maketitle


\section{Introduction}
\label{sec:intro}


\subsection{Notation}
\label{subsec:notation}

If $X$ is a variety endowed with a self-map $\varphi$, then the \emph{forward orbit} of a point $x\in X$ under the action of $\varphi$ is the set $\{\varphi^n(x):n\in\N\}$. We say that $x$ is \emph{preperiodic} under the action of $\varphi$ if its forward orbit is finite, i.e., there exist integers $0\leq m<n$ such that $\varphi^m(x)=\varphi^n(x)$ (where $\varphi^0$ is the identity map).  If one may take $m=0$, then we say that $x$ is \emph{periodic}.

Let $g(z)\in\Qbar[z]$ be a polynomial of degree $d\geq 3$. We say that $g$ is in \emph{normal form} (or equivalently, that it is \emph{normalized}) if it is monic and the coefficient of $z^{d-1}$ equals $0$. Every polynomial over $\Qbar$ is \emph{linearly conjugate} to a polynomial in normal form (i.e., there exists a linear polynomial $\mu(z)\in\Qbar[z]$ such that $\mu^{-1}\circ g\circ \mu$ is normalized); therefore, in this paper, we focus on normalized polynomials. Each normalized polynomial  of degree $d$ is represented  by a point in $\A^{d-1}$, as follows:
\begin{equation}
\label{eq:generic-family}
f_{\bt}(z)=z^d+t_1z^{d-2}+t_2z^{d-3}+\cdots+t_{d-2}z+t_{d-1},
\qquad
\bt=(t_1,\ldots,t_{d-1})\in\A^{d-1}.
\end{equation}
We fix (pairwise) distinct points $\alpha_1,\ldots,\alpha_d\in\Qbar$ and a target point $\beta\in\Qbar$. We are interested in the set of parameters $\bl$ for which $\beta$ lies in the forward orbit of each of the points $\alpha_i$ under the action of $f_{\bl}$. More precisely, we define
\begin{equation}
\label{eq:intro-collision-set}
C_f(\alpha_1,\ldots,\alpha_d;\beta)
:=
\left\{
\bl\in\A^{d-1}(\Qbar):
\begin{array}{c}
\text{for each }i\in\{1,\ldots,d\}\text{ there exists }n_i\in\N\\
\text{such that }f_{\bl}^{n_i}(\alpha_i)=\beta
\end{array}
\right\}.
\end{equation}
We stress that the integers $n_i$ are allowed to depend both on the starting point and on the parameter $\bl$. 


\subsection{Unlikely intersections in arithmetic dynamics}
\label{subsec:unlikely}

Given the action of a self-map $\varphi$ on an algebraic variety $X$, one may consider an \emph{unlikely event} $\mathcal E$ and ask what additional structure on the dynamical system $(X,\varphi)$ is forced if the event $\mathcal E$ occurs too often. This principle has led to several important questions in arithmetic dynamics. We mention two classical examples:
\begin{itemize}
\item given a subvariety $Y\subset X$, one asks whether the existence of a Zariski-dense set of preperiodic points of $\varphi$ on $Y$ forces $Y$ itself to be  preperiodic under another map dynamically related to $\varphi$; this is the setting of the dynamical Manin--Mumford problem (see \cite{Zha98,GTZ11,GT21});
\item given a point $x\in X$ and a subvariety $Y\subset X$, one asks whether an infinite intersection between the forward orbit of $x$ and $Y$ forces $Y$ to contain a positive-dimensional $\varphi$-periodic subvariety; this is the dynamical Mordell--Lang problem (see \cite{BGT16}).
\end{itemize}
For a general discussion of unlikely intersections in arithmetic geometry, we refer the reader to the excellent book of Zannier \cite{Zan12}; for surveys emphasizing the dynamical setting, see \cite{BIJMST19, Ghi17}. The setting of algebraic families gives rise to another natural unlikely-intersection problem. Given a family of maps together with several marked points, one may ask for parameters at which the marked points are simultaneously preperiodic. For one-parameter polynomial families this problem was initiated by Baker--DeMarco \cite{BD11}, with its arithmetic-geometric origins in the work of Masser--Zannier on simultaneous torsion in families of elliptic curves \cite{MZ10,MZ12}. Variants of this problem were subsequently studied in several settings; see, among others, \cite{BD13,GHT13,GHT15,FG22}.

It is equally natural to replace simultaneous preperiodicity by a prescribed common target. Thus, given a family $f_{\bt}$, starting points $\alpha_1,\ldots,\alpha_r$, and a target $\beta$, one asks for parameters $\bl$ such that
\begin{equation}
\label{eq:intro-collision-event}
f_{\bl}^{n_i}(\alpha_i)=\beta
\qquad\text{for each }i
\end{equation}
for suitable positive integers $n_i$. We refer to this as a \emph{collision of orbits} at the target $\beta$.

In a one-parameter family, the natural question is whether there are infinitely many collision parameters. In a higher-dimensional parameter space, however, infinitude is no longer the correct measure of largeness, since an infinite set may be contained in a proper algebraic subvariety of the parameter space (see Section~\ref{sec:general-family}). Accordingly, in the present paper we ask whether the collision set can be Zariski dense.


\subsection{Our main result}
\label{subsec:main}

We prove the following statement.

\begin{theorem}
\label{thm:main}
Let $d\geq 3$, and let $f_{\bt}$ be the generic family of monic normalized polynomials of degree $d$ from \eqref{eq:generic-family}. Let $\alpha_1,\ldots,\alpha_d\in\Qbar$ be pairwise distinct and let $\beta\in\Qbar$. Then $C_f(\alpha_1,\ldots,\alpha_d;\beta)$ is not Zariski dense in $\A^{d-1}$.
\end{theorem}

\begin{remark}
\label{rem:d=2}
We note that if $d=2$ (i.e., we have a $1$-parameter family of polynomials $f_t(z)=z^2+t$), then the corresponding conclusion from Theorem~\ref{thm:main} does not hold. Indeed, for any $c\in\Qbar^\ast$ and for any $b\in\Qbar$, the collision set $C_f(c,-c;b)$ is infinite. This is the case because both $c$ and $-c$ are mapped (after one iteration of our polynomial family) to the same point and then a simple application of \cite[Proposition~4.0.2]{GS26} yields the existence of infinitely many $\l\in\Qbar$ with the property that $f_\l^n(c^2+\l)=b$ for some suitable positive integer $n$. 
\end{remark}

Theorem~\ref{thm:main} is a special case of a more general result, which we state in Section~\ref{sec:general-family}. More precisely, for integers $d>k\geq 2$, we consider the $k$-parameter family
\begin{equation}
\label{eq:intro-general-family}
f_{\bt}(z)=z^d+t_1z^{k-1}+t_2z^{k-2}+\cdots+t_{k-1}z+t_k
\end{equation}
and prove that the collision set corresponding to $k+1$ pairwise distinct constant starting points and one constant target is not Zariski dense in $\A^k$. Taking $k=d-1$ in this more general result gives Theorem~\ref{thm:main}.

There is a simple dimension heuristic behind our theorem. Fix positive integers $n_1,\ldots,n_d$. Then the equations
$
f_{\bl}^{n_i}(\alpha_i)=\beta, \, i=1,\ldots,d, $
give $d$ algebraic conditions on the $(d-1)$-dimensional parameter space. Thus one expects the simultaneous collision of all $d$ marked orbits to be an unlikely event. The point of Theorem~\ref{thm:main} is that the same non-density conclusion remains true even though the \emph{hitting times} $n_1,\dots, n_d$ vary with the parameter $\bl$.

We also show that the precise placement of the parameters in \eqref{eq:intro-general-family} is relevant. If the parameters are attached to an arbitrary collection of monomials, then the resulting family may admit a nontrivial linear symmetry, and such a symmetry can force a Zariski-dense collision set; we give explicit examples in Section~\ref{sec:general-family} to show this phenomenon.


\subsection{Further connections between our theorem and other studied questions}
\label{subsec:previous}

The collision of orbits problem has already been studied in several arithmetic-dynamical settings. In our previous work \cite{GS26}, we considered arbitrary one-parameter families of normalized polynomials over number fields, with two starting points and one target, and obtained precise conditions under which the corresponding collision set is infinite. We studied an analogous problem for elliptic surfaces in \cite{GSell}. Related collision questions were considered for Drinfeld modules in \cite{Ghi24} and for polynomial families over fields of positive characteristic in \cite{AG26}.

Also, we mention that at its heart, the dimension-$1$ collision of orbits problem (i.e., for a $1$-parameter family of maps and two starting points) is a generalization of the classical \emph{GCD-problem}, which asks for understanding the variation of $\gcd\left(a^m-c, b^n-c\right)$ (for fixed positive integers $a,b,c$ and varying $m,n\in\N$). This GCD-problem originated in \cite{AR04, BCZ03} and variants of it were studied in diverse settings, from arbitrary $S$-units in number fields \cite{Luc05} to positive characteristic rings \cite{CZ13}. Also, motivated by a question from \cite{Sil04}, various geometric versions of the GCD-problem were considered in \cite{BCT24, GHT18}; furthermore, variants of this problem were considered in a purely dynamical setting (see \cite{HT17}), including in positive characteristic (see \cite{GHT17}).


\subsection{Strategy of proof and outline for our paper}
\label{subsec:plan}

We start in Section~\ref{sec:general-family} by stating the more general version of Theorem~\ref{thm:main}, for the family \eqref{eq:intro-general-family}, i.e., 
\begin{equation}
\label{eq:the_family}
f_\bt(z)=z^d + t_1z^{k-1}+t_2z^{k-2}+\cdots+t_{k-1}z+t_k\text{ (where $d>k\ge 2$).}
\end{equation}
We restrict to $k\ge 2$, since the non-density conclusion does not generally hold for one-parameter families (see Remark~\ref{rem:d=2}). The corresponding one-parameter collision problem is studied in \cite{GS26}. Also, in Section~\ref{sec:general-family}, we present examples showing that one cannot assign the parameters to arbitrary monomials in a normalized polynomial: root-of-unity symmetries may create Zariski-dense collision loci.

We continue in Section~\ref{sec:heights} with the height-theoretic setup. We first recall the Weil height and the canonical height associated to a specialization $f_{\bl}$ (for some $\bl\in\A^k(\Qbar)$). We then introduce the  map $\Phi_{\bl}$ on $\PP^k$ (which may be viewed as the diagonal extension of $f_\bl$ to $\PP^k$) and compare its canonical height with the canonical heights of the individual coordinatewise maps $f_\bl$. Applying \cite[Theorem~2.1]{GHT16}, we obtain
\begin{equation}
\label{eq:the_variation}
\hhat_{\Phi_{\bl}}(P)=\frac{1}{d}h(\bl)+O(1)
\end{equation}
for each constant $k$-tuple $P$ with pairwise distinct coordinates (where the implied constant from \eqref{eq:the_variation} is independent of $\bl\in\A^k(\Qbar)$). Let $\hhat_{f_{\bl}}$ denote the canonical height associated to a specialization $f_{\bl}$ (for some $\bl\in\A^{k}(\Qbar)$). If $f_{\bl}^{n_i}(\alpha_i)=\beta,$
then
\begin{equation}
\label{eq:intro-height-identity}
\hhat_{f_{\bl}}(\alpha_i)=d^{-n_i}\hhat_{f_{\bl}}(\beta).
\end{equation}
Thus, along a sequence of collision parameters for which the hitting times $n_i$ tend to infinity, one expects the starting points to have small canonical height, provided that the height of the target can be controlled. 

In Section~\ref{sec:small-collisions}, we assume that the collision set $C_f(\alpha_1,\dots,\alpha_{k+1};\beta)$ is Zariski dense. Using the identity \eqref{eq:intro-height-identity}, together with the height  variation estimate from Section~\ref{sec:heights} (see \eqref{eq:the_variation}), we show that the  height of the parameter $\bl$ is bounded (along some suitably chosen infinite sequence), the canonical height of the target must also be  bounded under $f_\bl$, and the canonical heights of all $k+1$ starting points must then tend to $0$.

Finally, in Section~\ref{sec:equidistribution}, we use the results of Section~\ref{sec:small-collisions} and show that our problem reduces to the problem studied in  \cite[\S5]{GHT16}. In particular, we explain that the equidistribution argument from \cite{GHT16} applies verbatim; this allows us to compare the height functions attached to two different $k$-tuples of marked points (among the $k+1$ starting  points $\alpha_1,\dots,\alpha_{k+1}$). This yields a rigidity statement saying that if any $k$ of the $k+1$ marked starting points are preperiodic at a given specialization $\bl\in\A^k(\Qbar)$, then so is the remaining marked starting point (with respect to the specialized map $f_\bl(z)$). We then use the non-density theorem for $k+1$ distinct points from \cite[Theorem~1.1]{GHN18}; the resulting contradiction proves our main theorem.


\section{More examples and the general form of our main result}
\label{sec:general-family}

We start by presenting an example showing that the set $C_f(\alpha_1,\dots, \alpha_{k+1};\beta)$ from  Theorem~\ref{thm:main} may be infinite. Indeed, consider the family of cubic polynomials
$$f_\bt(z)=z^3+t_1z+t_2,$$
along with the starting points $\alpha_1=0$, $\alpha_2=c$ and $\alpha_3=-c$ (for any given $c\in\Qbar^\ast$) and target point $\beta=0$. Then we consider the $1$-parameter family of polynomials 
$$g_t(z):= z^3+tz.$$
Using \cite[Proposition~4.0.2]{GS26}, we get that there exist infinitely many $\l\in\Qbar$ such that for some suitable $m\in\N$, we have that $g_\l^m(c)=0$ (note that the assumption that $c\ne 0$ yields that the hypotheses of \cite[Proposition~4.0.2]{GS26} are met). Hence, the set $C_f(0,c,-c;0)$ contains infinitely many points in $\bA^2(\Qbar)$ of the form $(\l,0)$.

We continue this section by stating the more general result  which will be studied throughout the remainder of our paper; Theorem~\ref{thm:main} is an immediate consequence of Theorem~\ref{thm:general}.


\subsection{The general family}
For integers $d>k\ge 2$, we let 
\begin{equation}\label{eq:family}
  f_{\bt}(z)
  =z^d+t_1z^{k-1}+t_2z^{k-2}+\cdots+t_{k-1}z+t_k,
  \qquad
  \bt=(t_1,\ldots,t_k)\in\A^k.
\end{equation}
For pairwise distinct $\alpha_1,\ldots,\alpha_{k+1}\in\Qbar$ and $\beta\in\Qbar$, set
\begin{equation}\label{eq:general-collision-set}
\begin{split}
  C_f(\alpha_1,\ldots,\alpha_{k+1};\beta)
  :=\{\bl\in\A^k(\Qbar):\;&\text{for every }i\in\{1,\ldots,k+1\}
  \text{ there exists }n_i\in\N\\
  &\text{such that }f_{\bl}^{n_i}(\alpha_i)=\beta\}.
\end{split}
\end{equation}

\begin{theorem}\label{thm:general}
Let $d>k\geq 2$, let $f_{\bt}$ be the family \eqref{eq:family}, let $\alpha_1,\ldots,\alpha_{k+1}\in\Qbar$ be pairwise distinct, and let $\beta\in\Qbar$. Then $C_f(\alpha_1,\ldots,\alpha_{k+1};\beta)$ is not Zariski dense in $\A^k$.
\end{theorem}

Taking $k=d-1$ in Theorem~\ref{thm:general} gives Theorem~\ref{thm:main}. 

\begin{remark}
\label{rem:d=2_again}
As explained in Remark~\ref{rem:d=2}, Theorem~\ref{thm:general} fails if $k=1$. Furthermore, the more general result from \cite[Theorem~1.3.1]{GS26} already provides a complete answer for the collision orbits problem of the $1$-parameter family of polynomials $f_t(z)=z^d+t$. Indeed, given the starting points $\alpha_1$ and $\alpha_2$ along with the target point $\beta$ (all three points living in $\Qbar$), \cite[Theorem~1.3.1]{GS26} yields that the set $C_f(\alpha_1,\alpha_2;\beta)$ of collision parameters is infinite if and only if $\alpha_1^d=\alpha_2^d$. The conclusions~(A)-(D) from \cite[Theorem~1.3.1]{GS26} (which guarantee the infinitude of $C_f(\alpha_1,\alpha_2;\beta)$) can only be met when $\alpha_1^d=\alpha_2^d$. This is because conditions~(A),~(C)~and~(D) never hold, while condition~(B) only occurs when $\alpha_1^d=\alpha_2^d$.
\end{remark}

Next, we explain why the precise form of the family \eqref{eq:family} is relevant. In particular, the restriction on the monomial support is important: if the parameters are attached to a different collection of monomials, a nontrivial linear symmetry may relate the marked orbits and force the collision set to be Zariski dense.

\subsection{A counterexample for more general families of polynomials}\label{subsec:counter-general}
One can easily construct generalizations of our next example; we prefer to give a concrete example of a $2$-parameter family of polynomials. So, we let
$$f_{\bt}(z)=z^6 + t_1z^3 + t_2,$$
along with starting points $\alpha_i:=c\zeta_3^i$ for $i=1,2,3$ (for some given $c\in\Qbar^\ast$, where $\zeta_3$ is a primitive third-root of unity) and target point $\beta:=b\in\Qbar$. It is immediate to see that 
$$f_\bt(\alpha_1)=f_\bt(\alpha_2)=f_\bt(\alpha_3)=c^6 + c^3\cdot t_1 + t_2.$$
For any given $\l_1\in\Qbar$, we consider the $1$-parameter family of polynomials
$$g^{(\l_1)}_t(z):=z^6 + \l_1 z^3 + t.$$
Then applying \cite[Proposition~4.0.2]{GS26} to starting (non-constant)  point $\alpha(t):=c^6 + c^3\l_1 + t$ and target point $\beta=b$, we obtain that there exist infinitely many $\l_2\in\Qbar$ such that $\beta$ is in the forward orbit of $\alpha(\l_2)$ under the action of $g^{(\l_1)}(z)$. In particular, this means that the original collision set of parameters $C_f(\alpha_1,\alpha_2,\alpha_3;\beta)\subset\A^2$ contains infinitely many points on each line $t_1=\l_1$ (for any given $\l_1\in\Qbar$). Therefore,  $C_f(\alpha_1,\alpha_2,\alpha_3;\beta)$ is indeed Zariski dense in $\A^2$.

\section{Heights and the general setup for our problem}
\label{sec:heights}


Throughout this section, let $d>k\geq 2$ be integers and we recall the notation for the family of polynomials from \eqref{eq:family}, i.e., $$f_\bt(z)=z^d+\sum_{i=1}^k t_i z^{k-i};$$
specializing each $t_i$ to $\l_i\in\Qbar$ yields a different map  $f_\bl(z)\in\Qbar[z]$  from our family of polynomials.

\subsection{Product formula fields}

We recall that a product formula field $K$  comes
equipped with a standard set $\Omega_K$ of absolute values $|\cdot|_v$ along with some positive integers $N_v$ (called \emph{local weights}) such that

\begin{equation}
\label{product formula equation 0}
\prod_{v\in \Omega_K} |x|^{N_v}_v = 1\quad\text{ for every $x\in K^*$.}
\end{equation}
In our paper, we will use two different types of product formula fields:
\begin{enumerate}
\item[(1)]  number fields (with the usual archimedean and nonarchimedean absolute values); and
\item[(2)] the rational function field $\Qbar(t_1,\dots, t_k)$ with the absolute values corresponding to the irreducible divisors of $\bP^k_{\Qbar}$ (after identifying $\Qbar(t_1,\dots, t_k)$ with the function field $\Qbar(\bP^k)$). 
\end{enumerate}
For more details, we refer the reader to \cite{BG06, lang}. Furthermore, when $K$ is a number field, then we fix an algebraic closure $\Qbar$ of $K$. Also, when $K$ is a number field, for each  $v \in \Omega_K$, we 
fix an extension of $|\cdot |_v$ to $\Qbar$.

\subsection{The Weil height for number fields}

Let $K$ be a number field and let $m\ge 1$ be an integer. Then for each finite extension $L$ of $K$ and for any point $P:=[x_m:\cdots :x_0]\in \bP^m(L)$, its Weil height is defined to be (see also \eqref{product formula equation 0} for the local weights $N_v$)
\begin{equation}
\label{eq:weil-height-projective}
h(P)=\frac{1}{[L:K]}\sum_{v\in\Omega_K}N_v\cdot \sum_{\substack{\sigma:L\into \Kbar\\ \sigma|_K=\id}} \log\left(\max\{|\sigma(x_m)|_v, \cdots, |\sigma(x_0)|_v\}\right).
\end{equation}
The above definition of height does not depend on the particular choice of the field $L$ containing each $x_i$. We also use the notation $h((x_m,\dots, x_1)):=h([x_m:\cdots :x_1:1])$ to denote the height of the point $(x_m, \dots, x_1)$ in the affine space $\bA^m$ embedded in the usual way in $\bP^m$.

\subsection{Weil height for points in a rational function field}

In the special case of the function field $K_0=\Qbar(t_1,\dots, t_k)$, for a point $P=[x_m:\cdots :x_0]\in\bP^m(K_0)$ (for some integer $m\ge 1$), assuming each $x_i\in \Qbar[t_1,\dots, t_k]$ and moreover, the polynomials $x_i$ are coprime, then we can simply define the Weil height of $P$ as follows: 
\begin{equation}
\label{eq:Weil_function_fields}
h(P)=\max_{i=0}^m \deg(x_i)\text{, where $\deg(\cdot )$ is the total degree function on $\Qbar[t_1,\dots, t_k]$.} 
\end{equation}
Again, for more details, we refer the reader to \cite{BG06, lang}.

\subsection{Canonical heights associated to maps defined over number fields}

Let $m\ge 1$, and let $f:\PP^m\lra \PP^m$ be an endomorphism of degree
$d\ge 2$ defined over $\Qbar$.  In \cite{CS93}, Call and Silverman defined the \emph{global canonical height} $\hhat_f(x)$ for each $x\in \PP^m(\Qbar)$ as
\begin{equation}
\label{definition canonical height number field}
\hhat_f(x)=\lim_{n\to\infty} \frac{h(f^n(x))}{d^n}.
\end{equation}
A classical application of Northcott's Theorem yields that 
\begin{equation}
\label{eq:Northcott}
\text{$\alpha$ is preperiodic for $f$ if and only if $\hhat_f(\alpha)=0$.}
\end{equation}

\subsection{Canonical height over function fields}
Given an endomorphism $f:\PP^m\lra \PP^m$ of degree 
$d\ge 2$ defined over the function field $K_0=\Qbar(t_1,\dots, t_k)$, we define similarly the canonical height of a point $\alpha\in \PP^m(K_0)$ as
\begin{equation}
\label{definition canonical height function field}
\hhat_f(\alpha)=\lim_{n\to\infty} \frac{h(f^n(\alpha))}{d^n}.
\end{equation}

In particular, regardless whether we work over a number field, or over a function field, if $f^n(\alpha)=\beta$, then
\begin{equation}\label{eq:collision-height}
  \hhat_f(\alpha)=d^{-n}\hhat_f(\beta).
\end{equation}

\subsection{The product map}
\label{subsec:product-variation}
Using the notation from Theorem~\ref{thm:general}, we let $\alpha_1,\dots, \alpha_k$ and $\beta$ in $\Qbar$; in particular, $\alpha_i\ne\alpha_j$ for $i\ne j$. We let $K$ be a fixed number field containing $\alpha_1,\dots, \alpha_k,\beta$. Corresponding to the family of polynomials $f_\bt$  as in \eqref{eq:family}, we let
\[
  F_{\bt}(X,X_0)
  =X^d+t_1X^{k-1}X_0^{d-k+1}+\cdots+t_{k-1}XX_0^{d-1}+t_kX_0^d
\]
be the homogenization of $f_{\bt}$. We associate to $f_{\bt}$ the coordinatewise endomorphism
\begin{equation}\label{eq:Phi}
\begin{split}
  \Phi_{\bt}:\PP^k&\longrightarrow\PP^k,\\
  [X_1:\cdots:X_k:X_0]
  &\longmapsto
  [F_{\bt}(X_1,X_0):\cdots:F_{\bt}(X_k,X_0):X_0^d].
\end{split}
\end{equation}
Once again, specializing each $t_i$ to some $\l_i\in\Qbar$ yields the corresponding morphism $F_\bl(X_1,X_0)$. 
The coordinate forms defining $\Phi_{\bl}$ (for any $\l_1,\dots, \l_k\in\Qbar$) have no common projective zero, so $\Phi_{\bl}$ (and also $\Phi_{\bt}$) is an endomorphism of degree $d$. On the affine chart $X_0=1$, we have 
\begin{equation}\label{eq:Phi-iterate}
  \Phi_{\bt}^n([x_1:\cdots:x_k:1])
  =[f_{\bt}^n(x_1):\cdots:f_{\bt}^n(x_k):1].
\end{equation}
For each $\bl\in\A^k(\Qbar)$, we write $\hhat_{\Phi_{\bl}}$ for the canonical height associated to $\Phi_{\bl}$.

\begin{lemma}\label{lem:product-comparison}
Let $K$ be a number field and let $\alpha_1,\ldots,\alpha_k\in K$. Then  for every $\bl=(\l_1,\dots, \l_k)\in\A^k(\Qbar)$, we have:
\begin{equation}\label{eq:product-comparison}
  \max_{1\leq i\leq k}\hhat_{f_{\bl}}(\alpha_i)
  \leq
  \hhat_{\Phi_{\bl}}([\alpha_1:\cdots:\alpha_k:1])
  \leq
  \sum_{i=1}^k\hhat_{f_{\bl}}(\alpha_i).
\end{equation}
\end{lemma}

\begin{proof}
Fix $\bl\in\A^k(\Qbar)$, and choose a finite extension $L/K$ containing all coordinates of $\bl$. For every absolute value $|\cdot|_v$ on $L$  and for each $y_1,\dots,y_k\in L$, we have
\begin{equation}
\label{eq:v_L}
  \max_{1\leq i\leq k}\left(\max\{1,|y_i|_v\}\right)
  \leq
  \max\{1,|y_1|_v,\ldots,|y_k|_v\}
  \leq
  \prod_{i=1}^k\max\{1,|y_i|_v\}.
\end{equation}
Taking logarithms in \eqref{eq:v_L} and summing over the absolute values on $L$ (extending the absolute values from $\Omega_K$ with the local weights and the usual normalization as in \eqref{eq:weil-height-projective}) gives
\begin{equation}
\label{eq:one_inequality}
  \max_{1\le i\le k} h(y_i)
  \leq h([y_1:\cdots:y_k:1])
  \leq \sum_{i=1}^k h(y_i).
\end{equation}
Apply inequalities \eqref{eq:one_inequality} to $f_{\bl}^n(\alpha_1),\ldots,f_{\bl}^n(\alpha_k)$ (for $n\in\mathbb{N}$), divide by $d^n$, and let $n\to\infty$; this yields the desired conclusion in Lemma~\ref{lem:product-comparison}.
\end{proof}

In particular, using also \eqref{eq:Northcott}, we deduce the following  consequence.
\begin{corollary}
\label{cor:product-zero}
With the notation as in Lemma~\ref{lem:product-comparison}, we have:
\begin{equation}\label{eq:product-zero}
  \hhat_{\Phi_{\bl}}([\alpha_1:\cdots:\alpha_k:1])=0
  \quad\Longleftrightarrow\quad
  \alpha_1,\ldots,\alpha_k\text{ are all preperiodic for }f_{\bl}.
\end{equation}
\end{corollary}

We will also employ in our proof the following easy fact.
\begin{lemma}
\label{lem:easy_fact}
For every constant point $P:=[\alpha_1:\dots : \alpha_k:1]\in\bP^k(\Qbar)$, we have that $\hhat_{\Phi_\bt}(P)=\frac{1}{d}$ (where $\hhat_{\Phi_\bt}$ is the corresponding canonical height of $\Phi_{\bt}$ computed with respect to the function field $\Qbar(\bt)=\Qbar(t_1,\dots, t_k)$). 
\end{lemma}

\begin{proof}
First, for every constant $\alpha\in\Qbar$, 
\begin{equation}\label{eq:generic-degree-growth}
  \deg_{\bt}\bigl(f_{\bt}^n(\alpha)\bigr)=d^{n-1}
  \qquad(n\in\N).
\end{equation}
Indeed, this is clear for $n=1$. If the degree at the $n$-th iterate is $d^{n-1}$, then the term $(f_{\bt}^n(\alpha))^d$ in the next iterate has total degree $d^n$ in $\bt=(t_1,\dots, t_k)$, whereas every remaining terms have degree at most
$$1+(k-1)d^{n-1}<d^n\text{ (because $d>k$).}$$ 
This proves \eqref{eq:generic-degree-growth} by induction. Consequently, for $P=[\alpha_1:\cdots:\alpha_k:1]\in\PP^k(\Qbar)$, its Weil height (computed for the function field $\Qbar(t_1,\dots, t_k)$) equals $h(\Phi_{\bt}^n(P))=d^{n-1}$ for each $n\ge 1$; hence
\begin{equation}\label{eq:generic-height}
  \hhat_{\Phi_\bt}(P)
  =\lim_{n\to\infty}\frac{d^{n-1}}{d^n}
  =\frac{1}{d},
\end{equation}
as desired in the conclusion of Lemma~\ref{lem:easy_fact}.
\end{proof}

\subsection{Variation of the canonical height}

Using Lemma~\ref{lem:easy_fact}, we can now state the precise form of the higher-dimensional specialization theorem of Ghioca--Hsia--Tucker \cite[Theorem~2.1]{GHT16}, which is key to our proof of Theorem~\ref{thm:general}. We also note that similar results for $1$-dimensional families were obtained in 
\cite{GM13, Ing13, Sil83, Tat83}.

\begin{theorem}[Ghioca--Hsia--Tucker]\label{thm:variation}
Let $\alpha_1,\ldots,\alpha_k\in\Qbar$ be pairwise distinct and put $P:=[\alpha_1:\cdots:\alpha_k:1]$. Then there exists a constant $C_P>0$ (depending only on $P$) such that, for every $\bl\in\A^k(\Qbar)$, we have: 
\begin{equation}\label{eq:variation-bounds}
  \frac1d h(\bl)-C_P
  \leq \hhat_{\Phi_{\bl}}(P)
  \leq \frac1d h(\bl)+C_P.
\end{equation}
\end{theorem}

\section{Points of small height}\label{sec:small-collisions}

Throughout this section, we work with the notation as in Theorem~\ref{thm:general}. Furthermore, we
assume, towards a contradiction, that
\begin{equation}\label{eq:ZD-assumption}
  C_f(\alpha_1,\ldots,\alpha_{k+1};\beta)
  \text{ is Zariski dense in }\A^k.
\end{equation}
A sequence $\{\bl_j\}\subset\A^k(\Qbar)$ is called \emph{generic} if every proper Zariski-closed subset of $\A^k$ contains only finitely many terms of the sequence. 

In order to state the main result of this section (see Proposition~\ref{prop:small-collision-sequence}), we recall the classical notation $O(1)$; more precisely, given a sequence $\{b_n\}_{n\ge 1}$ of nonnegative real numbers, we write $b_n=O(1)$ to denote the existence of a constant $c$ such that $b_n\le c$ for each $n\ge 1$.

\begin{proposition}\label{prop:small-collision-sequence}
Under the assumption \eqref{eq:ZD-assumption}, there exist a generic sequence
\[
  \{\bl_j\}_{j\geq1}
  \subseteq C_f(\alpha_1,\ldots,\alpha_{k+1};\beta)
\]
and for each $i=1,\dots, k+1$, there exist positive integers $\left\{n_{i,j}\right\}_{j\ge 1}$ such that
\begin{equation}\label{eq:hitting}
  f_{\bl_j}^{n_{i,j}}(\alpha_i)=\beta\text{ for each $i=1,\dots, k+1$ and for each $j\ge 1$}
\end{equation}
and also,  
\begin{equation}\label{eq:hitting-to-infinity}
  \lim_{j\to\infty} n_{i,j}=\infty\text{ for each $i=1,\dots, k+1$.}
\end{equation}
Moreover, we must also have that 
\begin{equation}\label{eq:parameter-bounded}
  h(\bl_j)=O(1),
\end{equation}
\begin{equation}\label{eq:target-bounded}
  \hhat_{f_{\bl_j}}(\beta)=O(1),
\end{equation}
and for each $i=1,\dots, k+1$, we have 
\begin{equation}\label{eq:marked-small}
  \lim_{j\to\infty}\hhat_{f_{\bl_j}}(\alpha_i)=0.
\end{equation}
\end{proposition}

\begin{proof}
For $i\in\{1,\ldots,k+1\}$ and $n\in\N$, let $V_{i,n}:=\{\bl\in\A^k:f_{\bl}^n(\alpha_i)=\beta\}.$
Each $V_{i,n}$ is a proper Zariski-closed subset of $\A^k$. Indeed, to see that it is proper, we note that $\bt\mapsto f_{\bt}^n(\alpha_i)$ is a polynomial of total degree $d^{n-1}$ (see equation~\eqref{eq:generic-degree-growth}); hence,  $f_{\bt}^n(\alpha_i)-\beta$ is not the zero polynomial.

Now, for $N\in\N$, set
\[
  E_N:=\bigcup_{i=1}^{k+1}\ \bigcup_{1\leq n\leq N}V_{i,n};
\]
each $E_N$ is a proper closed subset of $\A^k$. Now, let  $\left\{W_m\right\}_{m\ge 1}$ be an enumeration of the proper Zariski closed subsets of $\A^k_{\Qbar}$.  
Using \eqref{eq:ZD-assumption}, for each $j\ge 1$, choose
\begin{equation}\label{eq:generic-choice}
  \bl_j\in
  C_f(\alpha_1,\ldots,\alpha_{k+1};\beta)
  \setminus(E_j\cup W_1\cup\cdots\cup W_j).
\end{equation}
The resulting sequence is generic. For each $i$ and $j$, let  $n_{i,j}\in\N$ satisfying \eqref{eq:hitting}. Since $\bl_j\notin E_j$, we have $n_{i,j}>j$, which proves \eqref{eq:hitting-to-infinity}.

We write (for the sake of simplifying our notation)
\[
  f_j:=f_{\bl_j},\qquad
  \Phi_j:=\Phi_{\bl_j},\qquad
  H_j:=h(\bl_j),\qquad
  B_j:=\hhat_{f_j}(\beta).
\]
By \eqref{eq:hitting} and the functoriality of canonical height,
\begin{equation}\label{eq:alpha-height}
  \hhat_{f_j}(\alpha_i)=d^{-n_{i,j}}B_j
  \qquad(i=1,\ldots,k+1).
\end{equation}

Consider $P=[\alpha_1:\cdots:\alpha_k:1].$
Since $\alpha_1,\ldots,\alpha_k$ are pairwise distinct, Theorem~\ref{thm:variation} applies. Put $\varepsilon_j:=\sum_{i=1}^k d^{-n_{i,j}}$; by \eqref{eq:hitting-to-infinity}, we have 
\begin{equation}
\label{eq:epsilon_j}
\lim_{j\to\infty}\varepsilon_j=0. 
\end{equation} 
Combining the lower bound in \eqref{eq:variation-bounds} with Lemma~\ref{lem:product-comparison} and \eqref{eq:alpha-height}, we obtain that for some positive constant $C_P$ (independent of $j$), we have
\begin{equation}\label{eq:first-height-inequality}
  \frac1dH_j-C_P
  \leq \hhat_{\Phi_j}(P)
  \leq \sum_{i=1}^k\hhat_{f_j}(\alpha_i)
  =\varepsilon_jB_j.
\end{equation}

Next, we control $B_j$ in terms of $H_j$. Since the $k+1$ marked points are pairwise distinct, at most one of them equals $\beta$. We may therefore choose $k-1$ of the marked points, say $\gamma_2,\ldots,\gamma_k$, such that $\beta,\gamma_2,\ldots,\gamma_k$ are pairwise distinct. Let $Q=[\beta:\gamma_2:\cdots:\gamma_k:1].$
Lemma~\ref{lem:product-comparison} and the upper bound in \eqref{eq:variation-bounds} give
\begin{equation}\label{eq:target-vs-parameter}
  B_j
  \leq \hhat_{\Phi_j}(Q)
  \leq \frac1dH_j+C_Q
\end{equation}
for a positive constant $C_Q$ independent of $j$. Substituting \eqref{eq:target-vs-parameter} into \eqref{eq:first-height-inequality} yields
\[
  \frac1dH_j-C_P
  \leq
  \varepsilon_j\left(\frac1dH_j+C_Q\right),
\]
and hence
\begin{equation}\label{eq:height-rearranged}
  \frac{1-\varepsilon_j}{d}H_j
  \leq C_P+\varepsilon_jC_Q.
\end{equation}
For all sufficiently large $j$, we have $\varepsilon_j<1/2$ (see \eqref{eq:epsilon_j}), so \eqref{eq:height-rearranged} bounds $H_j$ uniformly. Enlarging the bound to include the finitely many remaining indices proves \eqref{eq:parameter-bounded}. Equation \eqref{eq:target-vs-parameter} then gives \eqref{eq:target-bounded}, and finally \eqref{eq:alpha-height} together with \eqref{eq:hitting-to-infinity} gives \eqref{eq:marked-small}.
\end{proof}

\section{Proof of Theorem~\ref{thm:general}}
\label{sec:equidistribution}

We work with the notation and hypotheses as in  Theorem~\ref{thm:general}. Furthermore, we will employ the notation from Section~\ref{subsec:product-variation} for the family of morphisms $\Phi_\bl$ (see \eqref{eq:Phi}) which extend the action of $f_\bl$ to an endomorphism $\bP^k$. 

\subsection{A rigidity result coming from an equidistribution theorem for points of small height}\label{subsec:rigidity}

Let $\alpha_1,\ldots,\alpha_{k+1}\in\Qbar$ be pairwise distinct and set
\begin{align}
  P^{(1)}&=[\alpha_k:\alpha_{k-1}:\cdots:\alpha_1:1]\in \bP^k(\Qbar), \label{eq:P-one}\\
P^{(2)}&=[\alpha_{k+1}:\alpha_{k-1}:\cdots:\alpha_1:1]\in \bP^k(\Qbar).\label{eq:P-two}
\end{align}
We first state the version of the rigidity argument which will be needed in our proof.

\begin{proposition}\label{prop:small-sequence-rigidity}
Suppose there exists a generic sequence $\{\bl_j\}\subset\A^k(\Qbar)$ such that
\begin{equation}\label{eq:two-small-products}
 \lim_{j\to\infty} \hhat_{\Phi_{\bl_j}}(P^{(1)})=0,
  \qquad
  \lim_{j\to\infty}\hhat_{\Phi_{\bl_j}}(P^{(2)}) =0.
\end{equation}
Then, for every $\bl\in\A^k(\Qbar)$,
\begin{equation}\label{eq:height-equality-products}
  \hhat_{\Phi_{\bl}}(P^{(1)})
  =
  \hhat_{\Phi_{\bl}}(P^{(2)}).
\end{equation}
In particular, for each $\bl\in\A^k(\Qbar)$, we have: 
\begin{equation}\label{eq:zero-locus-equivalence}
  \hhat_{\Phi_{\bl}}(P^{(1)})=0
  \quad\Longleftrightarrow\quad
  \hhat_{\Phi_{\bl}}(P^{(2)})=0.
\end{equation}
\end{proposition}

Proposition~\ref{prop:small-sequence-rigidity} is a consequence of the powerful theorem of Yuan \cite{Yua08} regarding the  equidistribution of points of small height with respect to a metrized line bundle. Yuan's theorem is part of a broader circle of results on small points --- see, among others, \cite{Bil97,BR06,CL06,FR06,Zha98}. 

\begin{proof}[Proof of Proposition~\ref{prop:small-sequence-rigidity}.]
Since the argument is already contained in the proof of \cite[Theorem~5.1]{GHT16}, we only point out the small difference between our hypothesis from Proposition~\ref{prop:small-sequence-rigidity} and the hypothesis of the aforementioned result from \cite{GHT16}. Indeed, in \cite[Theorem~5.1]{GHT16}, the authors used the hypothesis that there exists a generic sequence of points $\{\bl_j\}_j$ in $\bP^k(\Qbar)$ such that 
\begin{equation}
\label{eq:stronger}
\hhat_{\Phi_{\bl_j}}(P^{(1)})=\hhat_{\Phi_{\bl_j}}(P^{(2)})=0.
\end{equation}   
However, all one needs is the hypothesis from \eqref{eq:two-small-products} instead of the stronger hypothesis~\eqref{eq:stronger}. Indeed, exactly as constructed in the proof of   \cite[Theorem~5.1]{GHT16}, we have the two  metrized line bundles $\overline{\mathcal{L}}^{(1)}$ and $\overline{\mathcal{L}}^{(2)}$ on $\mathcal{O}_{\bP^k}(1)$, which are uniform limits of two sequences of algebraic ad\'elic metrized line bundles. Furthermore, for each $\bl=(\l_1,\l_2,\dots,\l_k)\in\A^k(\Qbar)$ (viewed as $[\l_1:\l_2:\dots:\l_k:1]\in\bP^k(\Qbar)$) and for each $j=1,2$, we have that (see \cite[equation~(5.2.3)]{GHT16})
\begin{equation}
\label{connection heights}
\hhat_{\overline{\cL}^{(j)}}(\bl)= d\cdot \hhat_{\Phi_{\bl}}\left(P^{(j)}\right).
\end{equation}
Combining equation \eqref{connection heights} with hypothesis~\eqref{eq:two-small-products} allows us to apply \cite[Corollary~4.3]{GHT15} (which is stated for generic sequences of small points with respect to two metrized line bundles) in order to conclude that (see \cite[equation~(5.2.9)]{GHT16}) 
\begin{equation}
\label{eq:two-heights-equal}
\hhat_{\overline{\cL}^{(1)}}(\bl) = \hhat_{\overline{\cL}^{(2)}}(\bl)\text{ for each $\bl\in\bP^k(\Qbar)$.}
\end{equation}
Equations~\eqref{eq:two-heights-equal} and \eqref{connection heights} yield the equality of the two canonical heights of $P^{(j)}$ (for $j=1,2$)  from equation \eqref{eq:height-equality-products} (for each $\bl\in\A^k(\Qbar)$). In particular, we get that $\hhat_{\Phi_{\bl}}(P^{(1)})=0$ if and only if $\hhat_{\Phi_{\bl}}(P^{(2)})=0$, i.e., the conclusion from equation~\eqref{eq:zero-locus-equivalence} holds.

This concludes our proof of Proposition~\ref{prop:small-sequence-rigidity}.
\end{proof}

\begin{corollary}\label{cor:preperiodicity-rigidity}
Assume that $C_f(\alpha_1,\ldots,\alpha_{k+1};\beta)$ is Zariski dense. Then, for every $\bl\in\A^k(\Qbar)$, if any $k$ of the points $\alpha_1,\ldots,\alpha_{k+1}$ are preperiodic for $f_{\bl}$, then all $k+1$ points are preperiodic for $f_{\bl}$.
\end{corollary}

\begin{proof}
Equation~\eqref{eq:zero-locus-equivalence} (see also Corollary~\ref{cor:product-zero}) yields that if  $\alpha_1,\dots, \alpha_k$ are preperiodic for some $f_\bl(z)$, then also $\alpha_{k+1}$ is preperiodic for $f_\bl(z)$. Since the order of the points $\alpha_1,\dots, \alpha_{k+1}$ is immaterial in the definition of  $C_f(\alpha_1,\ldots,\alpha_{k+1};\beta)$, we obtain the desired conclusion in Corollary~\ref{cor:preperiodicity-rigidity}.
\end{proof}


\subsection{Conclusion for our proof of Theorem~\ref{thm:general}}

Working under the hypotheses from Theorem~\ref{thm:general}, we argue by contradiction and therefore, assume $C_f(\alpha_1,\dots, \alpha_{k+1};\beta)$ is Zariski dense in $\bA^k$. Proposition~\ref{prop:small-collision-sequence} yields the existence of a generic sequence $\{\bl_j\}_j$ in $\bA^k$ with the property that
\begin{equation}
\label{eq:tending_all}
\lim_{j\to\infty} \hhat_{f_{\bl_j}}(\alpha_i)=0\text{ for $i=1,\dots, k+1$.}
\end{equation}
Equation \eqref{eq:tending_all} yields that for the two points $P^{(1)}$ and $P^{(2)}$ defined as in \eqref{eq:P-one} and \eqref{eq:P-two}, we have
\begin{equation}\label{eq:two-small-products_2}
 \lim_{j\to\infty} \hhat_{\Phi_{\bl_j}}(P^{(1)})=0,
  \qquad
  \lim_{j\to\infty}\hhat_{\Phi_{\bl_j}}(P^{(2)}) =0.
\end{equation}
Equation \eqref{eq:two-small-products_2} 
allows us to apply Proposition~\ref{prop:small-sequence-rigidity} and Corollary~\ref{cor:preperiodicity-rigidity} to obtain that for each $\bl\in\bA^k(\Qbar)$, if $\alpha_j$ is preperiodic for $f_\bl(z)$ for each $j=1,\dots, k$, then also $\alpha_{k+1}$ is preperiodic for $f_\bl(z)$; thus \cite[Theorem~2.2]{GHN18} holds. Then arguing identically as in the proof of \cite[Theorem~1.1]{GHN18}, we obtain a contradiction. 

Indeed, the main ingredient in the proof from \cite[Theorem~1.1]{GHN18} is \cite[Theorem~2.2]{GHN18}, which is precisely our Corollary~\ref{cor:preperiodicity-rigidity}. We briefly sketch the steps in the argument from \cite{GHN18}. For each function $\sigma:\{1,\dots, k-1\}\lra \{1,\dots, k-1\}$ (not necessarily injective), one considers the line $L_\sigma\subset \bA^k$ in the parameter space along which the following equations hold:
\begin{equation}
\label{eq:line}
\text{ for each $i=1,\dots, k-1$, }f_\bl(\alpha_i)=\alpha_{\sigma(i)}\text{ for $\bl\in L_\sigma$.}
\end{equation}
Specializing $f_\bt$ along the line $L_\sigma$ yields a  $1$-parameter family of polynomials $g^{(\sigma)}_t(z)$ (see \cite[equation~(3.5)]{GHN18}). Corollary~\ref{cor:preperiodicity-rigidity} yields that for each $\lambda\in\Qbar$, we have that $\alpha_k$ is preperiodic under the action of $g^{(\sigma)}_\lambda(z)$ if and only if $\alpha_{k+1}$ is preperiodic under the action of $g^{(\sigma)}_\lambda(z)$. Then the rest of the argument from the proof of \cite[Theorem~1.1]{GHN18} proceeds verbatim to yield that $\alpha_j=\alpha_{k+1}$ for some $j=1,\dots, k$. This contradicts the fact that the starting points $\alpha_1,\dots, \alpha_{k+1}$ are distinct. Therefore, indeed, the set $C_f(\alpha_1,\dots, \alpha_{k+1};\beta)$ is not Zariski dense in $\bA^k$, which concludes our proof of Theorem~\ref{thm:general}.


\end{document}